# Adaptive Third Order Adams-Bashforth Time Stepping for Extended Boussinesq Equations

Sasan Tavakkol[1, *], Sangyoung Son[1, 2], Patrick Lynett[1]

1- Department of Civil and Environmental Engineering, University of Southern California, Los Angeles, California, USA

2- School of Civil, Environmental and Architectural Engineering, Korea University, Republic of Korea

**ABSTRACT**

We develop the third-order adaptive Adams-Bashforth time stepping and the second-order finite difference equation for variable time steps. We incorporate these schemes in the Celeris Advent software to discretize and solve the 2D extended Boussinesq equations. This software uses a hybrid finite volume – finite difference scheme and leverages the GPU to solve the equations faster than real-time while concurrently visualizing them. We simulate several benchmarks using the adaptive time stepping scheme of Celeris Advent and demonstrate the capability of the software in modeling wave-breaking, wave runup, irregular waves, and rip currents. The adaptive scheme significantly improves the robustness of the model while providing faster computational performance.

## 1. Introduction

Numerical simulations are now essential tools in understanding any coastal phenomena ranging from wave behavior in ports to designing recreational surfing activities. Among nearshore models, the Boussinesq-type models have become the most popular approximations of the Navier–Stokes equations for coastal engineering, thanks to their ability to represent the major forces and interactions while requiring significantly less computational power compared to any 3D model. Nevertheless, these models are still computationally more expensive than their counterpart, non-linear shallow-water equations. This higher computational demand limits their application in low budget engineering projects. However, recent advances in computer hardware and software have lowered the barrier to entry for using the Boussinesq-type models. Recently, Tavakkol and Lynett [33] introduced a GPU accelerated software to solve the extended Boussinesq equations [23], called Celeris Advent. This software effectively democratized the use

---

[*] Corresponding author: tavakkol@usc.edu, currently at Google Research, New York, New York, USA

of Boussinesq-type models by letting the users run faster than real-time simulations on a consumer-level laptop and within a user-friendly interactive environment. In the current study, we explain and validate the adaptive time stepping scheme that we developed and incorporated in Celeris Advent [32].

The past three decades have seen a significant effort from the coastal research community towards developing Boussinesq-type models. Peregrine [25] derived the "standard" Boussinesq equations by assuming that both nonlinearity and frequency dispersion are weak and are in the same order of magnitude, therefore retaining only the lowest orders of nonlinearity and frequency dispersion terms. Because of this assumption, the standard Boussinesq equations are not applicable to very shallow or deep water. In shallow water, nonlinearity becomes more important than frequency dispersion as the wave gets closer to the shore, violating the assumption of same order nonlinear and dispersion effects. In the deep-water condition, the frequency dispersion cannot be considered weak for any depth greater than one-fifth of the equivalent deep-water wavelength, further limiting the application of the standard Boussinesq equations.

The deep-water restriction of the standard Boussinesq equations is often too limiting for engineering applications, especially where the incident wave energy spectrum consists of many frequency components. Several modified forms of Boussinesq equations have been successfully developed to extend their applications to deeper water depth (e.g., [7, 22, 24]). Among these extended Boussinesq equations, those introduced by Madsen and Sørensen [23] and Nwogu [24] are widely in use.

The weak nonlinearity restriction which limits the applicability of modified Boussinesq equations in very shallow waters is removed by eliminating this assumption in fully non-linear models (e.g., [17, 43]). FUNWAVE [9] and COULWAVE [19] are widely-used numerical implementation of the fully non-linear Boussinesq equations. These models have proven themselves successful in a wide variety of applications such as wave runup [21], wave-current interaction [28], wave generation by underwater landslides [18], rip and longshore currents [5], etc.

Fully non-linear models are known to better represent steep waves in shallow water and are shown to agree better with controlled laboratory experiments as well as with analytical solutions; however, their application in hindcasting or forecasting a real-world field site may not benefit from their higher-order accuracy because of uncertainties in the field site conditions (e.g., in boundary condition and bathymetry). Considering that these models are also computationally more expensive than the weakly non-linear Boussinesq models, their application in real-world scenarios might not be always justifiable. Bearing these facts in mind, we chose to solve the extended Boussinesq equations introduced by Madsen and Sørensen [23] in Celeris Advent and achieved sufficient accuracy with faster than real-time simulation speed.

Celeris Advent was originally developed to use a fixed third-order Adams-Bashforth time stepping scheme as the predictor step and an optional fixed fourth-order Adams-Moulton time stepping scheme as the corrector step. We later dropped the correction step as the predictor step proved to be sufficiently accurate, provided that an adequately small time step is chosen. Celeris

Advent with the fixed time step is validated for wave and current simulation in [33, 36] by the authors of the software as well as in [3, 26, 27] by other researchers. In this paper we describe the development of our new numerical time stepping scheme which accepts variable time step values and therefore allows the model to keep the maximum local CFL number constant by adaptively calculating the required time step value over the time. We developed third-order adaptive equations for the Adams-Bashforth time stepping scheme and incorporated them in Celeris Advent [32]. The need for an adaptive time stepping formulation arose from our observations of model instability in experiments with runup on steep surfaces. In these cases, the flow velocity grows rapidly and the subsequent increase in the local Froude number leads to instability. We validate the new adaptive scheme incorporated in Celeris Advent (v.1.3.4). The open source code and its compiled version are available to download at www.celeria.org.

Shi et al. [29] introduced a high-order adaptive time stepping solver for Boussinesq-type equations using Runge–Kutta time stepping. This scheme utilizes a fixed time step throughout the multi-level time stepping phase and adjusts the time step only for the next time level according to the CFL number. Therefore, this scheme is not theoretically fully adaptive. To the best of our knowledge, we are the first to introduce a third order adaptive time stepping scheme to solve Boussinesq equations [32]. We developed our third order adaptive scheme for Adams-Bashforth time stepping in a general format such that it can be used to solve equations other than the Boussinesq equations as well.

This paper is organized as follows. We describe the extended Boussinesq equations and a specific rearrangement in Section 2. In Section 3, we briefly explain our hybrid finite volume – finite difference scheme to solve the Boussinesq-type equations in space but explain in more detail their solution on an adaptive time grid. We give a short explanation on the development of Celeris Advent in Section 4, as the details are given previously in [32, 33]. In Section 5 we demonstrate the capability of our model by applying it on four benchmarks including wave breaking, wave runup, irregular waves, and rip currents. Conclusion, acknowledgments, and references make up the last three sections of this paper.

## 2. Extended Boussinesq equations

Celeris solves the extended Boussinesq equations derived by Madsen and Sørensen [23]. These equations for 2DH flow read as

$$\begin{bmatrix} \eta \\ P \\ Q \end{bmatrix}_t + \begin{bmatrix} P \\ P^2/h \\ PQ/h \end{bmatrix}_x + \begin{bmatrix} Q \\ PQ/h \\ Q^2/h \end{bmatrix}_y + \begin{bmatrix} 0 \\ gh\eta_x + \psi_1 + f_1 \\ gh\eta_y + \psi_2 + f_2 \end{bmatrix} = 0 \qquad (1)$$

where $\eta$ is the free surface elevation measured from the still water surface elevation, $h$ is the total water depth, $P$ and $Q$ are the depth-integrated mass fluxes in $x$ and $y$ directions, respectively, $g$ is the gravitational acceleration coefficient, and $f_1$ and $f_2$ are the bottom friction terms. Subscripts $x$ and $y$ denote spatial differentiation, with respect to the corresponding direction, and subscript $t$ denotes temporal differentiation. Finally, $\psi_1$ and $\psi_2$ are the modified dispersive terms defined as

$$\psi_1 = -\left(B + \frac{1}{3}\right)d^2(P_{xxt} + Q_{xyt}) - Bgd^3(\eta_{xxx} + \eta_{xyy})$$
$$- dd_x\left(\frac{1}{3}P_{xt} + \frac{1}{6}Q_{yt} + 2Bgd\eta_{xx} + Bgd\eta_{yy}\right) \quad (2)$$
$$- dd_y\left(\frac{1}{6}Q_{xt} + Bgd\eta_{xy}\right)$$

$$\psi_2 = -\left(B + \frac{1}{3}\right)d^2(P_{xyt} + Q_{yyt}) - Bgd^3(\eta_{yyy} + \eta_{xxy})$$
$$- dd_y\left(\frac{1}{3}Q_{yt} + \frac{1}{6}P_{xt} + 2Bgd\eta_{yy} + Bgd\eta_{xx}\right) \quad (3)$$
$$- dd_x\left(\frac{1}{6}P_{yt} + Bgd\eta_{xy}\right)$$

where $d$ is the still water depth and $B$ is the calibration coefficient for dispersion properties of the equations. We use $B=1/15$ as suggested originally in [22] and widely adopted thereafter.

The modified dispersive terms ($\psi_1$ and $\psi_2$) approach to zero as $d$ decreases to zero. This is a favorable property because as a wave approaches the shore (i.e., the still water depth decreases), it gets steeper, driving the waveheight ($H$) to still water depth ratio $\varepsilon = H/d$ higher and the square still water depth to wavelength ($L$) ratio $\mu^2 = (d/L)^2$ lower. Note that $\varepsilon$ and $\mu^2$ are representatives of the nonlinear effects and the dispersive effects, respectively. This progressive change in the values, as a wave gets closer to the shore, invalidate the underlying assumption of $O(\varepsilon)=O(\mu^2)$ for the derivation of Boussinesq equations and pushes the governing equations to the range where NLSW equations suit better. For $d = 0$, the extended Boussinesq equations, Eq. (1), reduces to the NLSW equations. In areas where the still water surface elevation is not defined, such as on lands above the sea level, we set $d = 0$, so the solver automatically switches to NLSW equations. The extended Boussinesq equations provide sufficiently accurate linear dispersion and shoaling characteristics for values of $kd < 3$, where $k$ is the wavenumber.

We rewrite Eq. (1) in a conservative form which is suitable to the applications of finite volume method. Expressing the free surface elevation as $\eta=h-d$ we have

$$\eta_t = h_t - d_t; \quad \eta_x = h_x - d_x; \quad \eta_y = h_y - d_y. \quad (4)$$

We let $b$ denote the bottom elevation from a fixed datum and $w_s$, a constant number, denote the still water elevation from this datum. Therefore, we have $d=w_s-b$. Since $w_s$ is constant in space and time, the derivative of the still water depth and the bottom elevation becomes equal, but with a negative sign. Furthermore, the temporal derivative of the still water depth, $d_t$, becomes zero assuming a constant bottom elevation in time.

We make a variable change by introducing $w=h+b$, where $w$ is the water surface elevation, measured from the fixed datum. This variable change helps us employ well-balanced numerical schemes for discretization of the advective terms, which is discussed in detail in [16]. Using the new notations, we rewrite Eq. (1) as

$$\mathbf{U}_t + \mathbf{F}(\mathbf{U})_x + \mathbf{G}(\mathbf{U})_y + \mathbf{S}(\mathbf{U}) = 0 \quad (5)$$

where newly introduced variables are

$$\mathbf{U} = \begin{bmatrix} w \\ P \\ Q \end{bmatrix}$$

$$\mathbf{F(U)} = \begin{bmatrix} P \\ P^2/(w-b) + 1/2\, g(w-b)^2 \\ PQ/w - b \end{bmatrix}$$

$$\mathbf{G(U)} = \begin{bmatrix} Q \\ PQ/(w-b) \\ Q^2/(w-b) + 1/2\, g(w-b)^2 \end{bmatrix}$$

$$\mathbf{S(U)} = \begin{bmatrix} 0 \\ (w-b)b_x + \psi_1 + f_1 \\ (w-b)b_y + \psi_2 + f_2 \end{bmatrix}$$

In Eq. (5), $\mathbf{U}$ is the conservative variables vector, $\mathbf{F(U)}$ and $\mathbf{G(U)}$ are the advective flux vectors, and $\mathbf{S(U)}$ is the source term which includes bottom slope, friction, and dispersive terms.

## 3. Numerical schemes

We use a hybrid finite volume – finite difference scheme on a uniform spatial Cartesian grid which we introduced in [33] and refer to it as TL17. We developed TL17 following similar works in [8] and [41]. In this scheme, the NLSW subset of the extended Boussinesq equations, Eq. (1), is discretized using a second-order well-balanced positivity preserving central-upwind scheme introduced by Kurganov and Petrova [16]. This scheme, known as KP07, is a finite volume method (FVM) to solve the Saint-Venant system of shallow water equations. The modified dispersive terms are discretized using the central finite difference method (FDM).

### 3.1 Spatial discretization

Following Wei and Kirby [42], Eq. (5) can be rearranged as

$$w_t = E(P, Q) \tag{6}$$
$$U_t^* = F(h, P, Q) + [F^*(Q)]_t \tag{7}$$
$$V_t^* = G(h, P, Q) + [G^*(P)]_t \tag{8}$$

where

$$U^* = P - \frac{1}{3}dd_x P_x - \left(B + \frac{1}{3}\right)d^2 P_{xx} \tag{9}$$

$$V^* = Q - \frac{1}{3}dd_y Q_y - \left(B + \frac{1}{3}\right)d^2 Q_{yy} \tag{10}$$

$$E(P, Q) = -(P_x + Q_y) \tag{11}$$

$$F(w, P, Q) = -\left(\frac{P^2}{w-b} + \frac{g(w-b)^2}{2}\right)_x - \left(\frac{PQ}{w-b}\right)_y - g(w-b)b_x - f_1 \\ + Bgd^3(\eta_{xxx} + \eta_{xyy}) + Bgd^2(d_x(2\eta_{xx} + \eta_{yy}) + d_y\eta_{xy}) \tag{12}$$

$$G(w, P, Q) = -\left(\frac{Q^2}{w-b} + \frac{g(w-b)^2}{2}\right)_y - \left(\frac{PQ}{w-b}\right)_x - g(w-b)b_y - f_2 \qquad (13)$$
$$+ Bgd^3(\eta_{yyy} + \eta_{xxy}) + Bgd^2(d_y(2\eta_{yy} + \eta_{xx}) + d_x\eta_{xy})$$

$$F^*(Q) = \frac{1}{6}dd_xQ_y + \frac{1}{6}dd_yQ_x + \left(B + \frac{1}{3}\right)d^2Q_{xy} \qquad (14)$$

$$G^*(Q) = \frac{1}{6}dd_xP_y + \frac{1}{6}dd_yP_x + \left(B + \frac{1}{3}\right)d^2P_{xy} \qquad (15)$$

The left-hand side terms in Eq. (6)-(8) are discretized in time, $[F^*(Q)]_t$ and $[G^*(P)]_t$ are evaluated by extrapolation in time, and the rest of the terms on the right hand side are known in the current time step. This rearrangement allows us to rewrite Eq. (5) as ODE's in time.

As mentioned before, we use KP07 to solve the NLSW subset of the extended Boussinesq equations. We chose this scheme because it is well-balanced (i.e., preserves stationary steady states) and guarantees the positivity of the computed fluid depth. Furthermore, it naturally supports a dry state, with no need to keep track of the wet-dry front, and it can accommodate discontinuous bottom topography. These qualities were required to develop Celeris Advent as an interactive solver. TL17 uses KP07 as the FVM solver and adds the dispersive terms discretized by central FDM to the source term in KP07. The spatial domain is discretized by rectangular cells with fixed sizes of $\Delta x$ and $\Delta y$. Each cell plays the role of a control volume for the FVM discretization. Cell centers and their corresponding cell averaged values are used as the grid points in FDM. The details of the KP07 and TL17 solver are explained in [16] and [33], respectively.

## 3.2 Time Integration

### 3.2.1 Uniform Time stepping

Uniform time integration is performed by the third-order Adams-Bashforth scheme which reads as

$$w_{ij}^{n+1} = w_{ij}^n + \frac{\Delta t}{12}\left(23E_{ij}^n - 16E_{ij}^{n-1} + 5E_{ij}^{n-2}\right) \qquad (16)$$

$$U_{ij}^{*n+1} = U_{ij}^{*n} + \frac{\Delta t}{12}\left(23F_{ij}^n - 16F_{ij}^{n-1} + 5F_{ij}^{n-2}\right) + 2F_{ij}^{*n} - 3F_{ij}^{*n-1} + F_{ij}^{*n-2} \qquad (17)$$

$$V_{ij}^{*n+1} = V_{ij}^{*n} + \frac{\Delta t}{12}\left(23G_{ij}^n - 16G_{ij}^{n-1} + 5G_{ij}^{n-2}\right) + 2G_{ij}^{*n} - 3G_{ij}^{*n-1} + G_{ij}^{*n-2} \qquad (18)$$

where the superscripts denote the step number in time, with *n* being the last step with known values. This time stepping is explicit in time, meaning that all the variables on the right-hand side of the equations are known. Since the variables at previous time steps are not defined in the very first two time steps of the simulation (i.e., *n*=1 and *n*=2), a first order Euler time integration is used to bootstrap the simulation until *n*=3.

The water surface elevation, $w^{n+1}$, is directly calculated from Eq. (16). However, to calculate the flux terms, $P^{n+1}$ and $Q^{n+1}$ the following set of implicit equations need to be solved:

$$A_{ij}^x P_{i-1,j} + B_{ij}^x P_{ij} + C_{i,j}^x P_{i+1,j} = U_{ij}^* \tag{19}$$

$$A_{ij}^y Q_{i,j-1} + B_{ij}^y Q_{ij} + C_{ij}^y Q_{i,j+1} = V_{ij}^* \tag{20}$$

where

$$\begin{aligned} A^\alpha &= \frac{dd_\alpha}{6\Delta\alpha} - \left(B + \frac{1}{3}\right)\frac{d^2}{\Delta\alpha^2}, \quad B^\alpha = 1 + 2\left(B + \frac{1}{3}\right)\frac{d^2}{\Delta\alpha^2}, \quad C^\alpha \\ &= -\frac{dd_\alpha}{6\Delta\alpha} - \left(B + \frac{1}{3}\right)\frac{d^2}{\Delta\alpha^2} \end{aligned} \tag{21}$$

The coefficient matrices in Eq. (19) and Eq. (20) are of tridiagonal form. We adopted the Cyclic Reduction (CR) method to efficiently solve these set of equations on the GPU in Celeris Advent.

### 3.2.2 Adaptive Time stepping

In the adaptive mode, the software keeps the maximum local CFL at a constant value, by using a variable time step. We define the CFL number for Celeris Advent as

$$\text{CFL} = \Delta t \times \underset{ij}{\text{MAX}}\left(\frac{\underset{ij}{\text{MAX}}(|u_{ij} \pm c_{ij}|)}{\Delta x}, \frac{\underset{ij}{\text{MAX}}(|v_{ij} \pm c_{ij}|)}{\Delta y}\right), \quad c_{ij} = \sqrt{gh_{ij}} \tag{22}$$

where $c$ is the wave celerity in shallow water. The theoretical stability condition for KP07, and thus for TL17, is CFL < 0.25 [16], however, in practice we often use a 0.5 safety factor and keep CFL smaller than 0.125. In the adaptive time stepping, Celeris calculates the next time step size, $\Delta t$, from Eq. (22) for a given constant CFL number and using the velocity and celerity values of the current time step.

#### 3.2.2.1 Third-order Adaptive Adams-Bashforth Equation

We aim to solve the following ODE:

$$X_t = f(t, X), \quad X(t_0) = X^0 \tag{23}$$

where $X_t$ denotes derivative of $X$ with respect to $t$, and $t_0$ denotes $t = 0$. Let $X^{i+1}$ denote our target variable at the next time step, $t_{i+1}$, and $X^i$ denote the same variable in the current time step, $t_i$. We can make an approximation of $f(t, X)$ by the third-degree polynomial, $p(t)$, such that:

$$p(t_{i-s}) = f(t_{i-s}, X^{i-s}), \quad for\ s = 0, 1, and\ 2 \tag{24}$$

Employing the Lagrange formula for polynomial interpolation we have:

$$p(t) = \sum_{j=i-2}^{i} \left( \prod_{\substack{k=i-2 \\ k \neq j}}^{i} \frac{t - t_k}{t_j - t_k} \right) X_t^j \qquad (25)$$

where $X_t^j = f(t_j, X^j)$. Now we can write:

$$X^{i+1} = X^i + \int_{t_i}^{t_{i+1}} p(t)\, dt \qquad (26)$$

$$X^{i+1} = X^i + \int_{t_i}^{t_{i+1}} \left( \sum_{j=i-2}^{i} \left( \prod_{\substack{k=i-2 \\ k \neq j}}^{i} \frac{t - t_k}{t_j - t_k} \right) X_t^j \right) dt \qquad (27)$$

$$\begin{aligned} X^{i+1} = X^i &+ \left( \int_{t_i}^{t_{i+1}} \frac{t - t_{i-1}}{t_{i-2} - t_{i-1}} \times \frac{t - t_i}{t_{i-2} - t_i}\, dt \right) X_t^{i-2} \\ &+ \left( \int_{t_i}^{t_{i+1}} \frac{t - t_{i-2}}{t_{i-1} - t_{i-2}} \times \frac{t - t_i}{t_{i-1} - t_i}\, dt \right) X_t^{i-1} \\ &+ \left( \int_{t_i}^{t_{i+1}} \frac{t - t_{i-2}}{t_i - t_{i-2}} \times \frac{t - t_{i-1}}{t_i - t_{i-1}}\, dt \right) X_t^i \end{aligned} \qquad (28)$$

Using the sliding technique to substitute $t$ with $t+t_i$ and introducing $\Delta t_i = t_{i+1} - t_i$, $\Delta t_{i-1} = t_i - t_{i-1}$, and $\Delta t_{i-2} = t_{i-1} - t_{i-2}$ we can write:

$$\begin{aligned} X^{i+1} = X^i &+ \left( \int_{0}^{\Delta t_i} \frac{t + \Delta t_{i-1}}{\Delta t_{i-2}} \times \frac{t}{\Delta t_{i-1} + \Delta t_{i-2}}\, dt \right) X_t^{i-2} \\ &- \left( \int_{0}^{\Delta t_i} \frac{t + (\Delta t_{i-1} + \Delta t_{i-2})}{\Delta t_{i-2}} \times \frac{t}{\Delta t_{i-1}}\, dt \right) X_t^{i-1} \\ &+ \left( \int_{0}^{\Delta t_i} \frac{t + (\Delta t_{i-1} + \Delta t_{i-2})}{\Delta t_{i-1} + \Delta t_{i-2}} \times \frac{t + \Delta t_{i-1}}{\Delta t_{i-1}}\, dt \right) X_t^i \end{aligned} \qquad (29)$$

Finally, after integration we have:

$$X^{i+1} = X^i + \frac{\Delta t_i}{6}\left[\left(\frac{\Delta t_i}{\Delta t_{i-1}} \times \frac{2\Delta t_i + 6\Delta t_{i-1} + 3\Delta t_{i-2}}{\Delta t_{i-1} + \Delta t_{i-2}} + 6\right)X_t^i \right.$$
$$\left. - \left(\frac{\Delta t_i}{\Delta t_{i-1}} \times \frac{2\Delta t_i + 3\Delta t_{i-1} + 3\Delta t_{i-2}}{\Delta t_{i-2}}\right)X_t^{i-1} \right. \tag{30}$$
$$\left. + \left(\frac{\Delta t_i}{\Delta t_{i-2}} \times \frac{2\Delta t_i + 3\Delta t_{i-1}}{\Delta t_{i-1} + \Delta t_{i-2}}\right)X_t^{i-2}\right]$$

Eq. (30) is the third order adaptive Adams-Bashforth time integration equation. As a correctness check, applying $\Delta t_i = \Delta t_{i-1} = \Delta t_{i-2} = \Delta t$ in this equation yields to the same third order Adams-Bashforth equation we use for uniform time stepping:

$$X^{i+1} = X^i + \frac{\Delta t}{12}[23 X_t^i - 16 X_t^{i-1} + 5 X_t^{i-2}]$$

### 3.2.2.2 Variable-step Second Order Finite Difference Equations

To solve Eq. (7) and Eq. (8) using Eq. (30) we also need to derive the second order finite difference discretization equation for variable time steps. Let's approximate $Y$ with $\hat{Y}$, and its derivative, $Y_t$, with $\hat{Y}_t$. We use a polynomial approximation for $Y$ such that it satisfies:

$$\hat{Y}(t_{i-s}) = Y(t_{i-s}), \quad for\ s = 0, 1, and\ 2 \tag{31}$$

Let's define:

$$\hat{Y} = c_0 + c_1(t - t_i) + c_2(t - t_i)^2 \tag{32}$$

The finite difference approximation of $Y_t^i$ yields:

$$Y_t^i \approx \hat{Y}_t^i = c_1 \tag{33}$$

To meet the conditions in Eq. (31), we must have:

$$\begin{bmatrix} 1 & 0 & 0 \\ 1 & -\Delta t_{i-1} & (\Delta t_{i-1})^2 \\ 1 & -(\Delta t_{i-1} + \Delta t_{i-2}) & (\Delta t_{i-1} + \Delta t_{i-2})^2 \end{bmatrix} \begin{bmatrix} c_0 \\ c_1 \\ c_2 \end{bmatrix} = \begin{bmatrix} Y^i \\ Y^{i-1} \\ Y^{i-2} \end{bmatrix} \tag{34}$$

Inverting the matrix of coefficients in (34), we can find the value of $c_1$:

$$Y_t^i \approx \hat{Y}_t^i = c_1 = \frac{2\Delta t_{i-1} + \Delta t_{i-2}}{\Delta t_{i-1}(\Delta t_{i-1} + \Delta t_{i-2})} Y^i - \frac{\Delta t_{i-1} + \Delta t_{i-2}}{\Delta t_{i-1}\Delta t_{i-2}} Y^{i-1}$$
$$+ \frac{\Delta t_{i-1}}{\Delta t_{i-2}(\Delta t_{i-1} + \Delta t_{i-2})} Y^{i-2} \tag{35}$$

For the special case where $\Delta t_{i-1} = \Delta t_{i-2} = \Delta t$

$$Y_t^i = \frac{3}{2\Delta t} Y^i - \frac{2}{\Delta t} Y^{i-1} + \frac{1}{2\Delta t} Y^{i-2} = \frac{1}{2\Delta t}(3Y^i - 4Y^{i-1} + Y^{i-2})$$

which is the well-known 2$^{nd}$ order backward difference.

Similarly, to derive $Y_t^{i-1}$, we can use the following polynomial approximation and solve it to find $\hat{Y}_t^{i-1}$:

$$\hat{Y} = c_0 + c_1(t - t_{i-1}) + c_2(t - t_{i-1})^2 \tag{36}$$

$$\begin{bmatrix} 1 & \Delta t_{i-1} & (\Delta t_{i-1})^2 \\ 1 & 0 & 0 \\ 1 & -\Delta t_{i-2} & (\Delta t_{i-2})^2 \end{bmatrix} \begin{bmatrix} c_0 \\ c_1 \\ c_2 \end{bmatrix} = \begin{bmatrix} Y^i \\ Y^{i-1} \\ Y^{i-2} \end{bmatrix} \tag{37}$$

$$Y_t^{i-1} \approx \frac{\Delta t_{i-2}}{\Delta t_{i-1}(\Delta t_{i-1} + \Delta t_{i-2})} Y^i + \frac{\Delta t_{i-1} - \Delta t_{i-2}}{\Delta t_{i-1}\Delta t_{i-2}} Y^{i-1} - \frac{\Delta t_{i-1}}{\Delta t_{i-2}(\Delta t_{i-1} + \Delta t_{i-2})} Y^{i-2} \tag{38}$$

For the special case where $\Delta t_{i-1} = \Delta t_{i-2} = \Delta t$

$$Y_t^{i-1} = \frac{1}{2\,\Delta t} Y^i + \frac{0}{\Delta t} Y^{i-1} + \frac{1}{2\Delta t} Y^{i-2} = \frac{1}{2\Delta t}(Y^i - Y^{i-2})$$

which is the well-known 2$^{nd}$ order central difference.

To derive $Y_t^{i-2}$ we can write:

$$\hat{Y} = c_0 + c_1(t - t_{i-2}) + c_2(t - t_{i-2})^2 \tag{39}$$

$$\begin{bmatrix} 1 & \Delta t_{i-1} + \Delta t_{i-2} & (\Delta t_{i-1} + \Delta t_{i-2})^2 \\ 1 & \Delta t_{i-2} & (\Delta t_{i-2})^2 \\ 1 & 0 & 0 \end{bmatrix} \begin{bmatrix} c_0 \\ c_1 \\ c_2 \end{bmatrix} = \begin{bmatrix} Y^i \\ Y^{i-1} \\ Y^{i-2} \end{bmatrix} \tag{40}$$

$$Y_t^{i-2} = -\frac{\Delta t_{i-2}}{\Delta t_{i-1}(\Delta t_{i-1} + \Delta t_{i-2})} Y^i + \frac{\Delta t_{i-1} + \Delta t_{i-2}}{\Delta t_{i-1}\Delta t_{i-2}} Y^{i-1} - \frac{\Delta t_{i-1} + 2\Delta t_{i-2}}{\Delta t_{i-2}(\Delta t_{i-1} + \Delta t_{i-2})} Y^{i-2} \tag{41}$$

For the special case where $\Delta t_{i-1} = \Delta t_{i-2} = \Delta t$

$$Y_t^{i-2} = -\frac{1}{2\,\Delta t} Y^i + \frac{2}{\Delta t} Y^{i-1} - \frac{3}{2\Delta t} Y^{i-2} = -\frac{1}{2\Delta t}(Y^i - 4Y^{i-1} + 3Y^{i-2})$$

which is the well-known 2$^{nd}$ order forward difference.

### 3.2.2.3 Prediction Equations

Now, we have the tools to derive the prediction equations of $w$, $U^*$, and $V^*$ with adaptive time stepping. Deriving the equation for $w$ is straightforward and yields:

$$w_{ij}^{n+1} = w_{ij}^n + \frac{\Delta t_n}{6}\left[\left(\frac{\Delta t_n}{\Delta t_{n-1}} \times \frac{2\Delta t_n + 6\Delta t_{n-1} + 3\Delta t_{n-2}}{\Delta t_{n-1} + \Delta t_{n-2}} + 6\right)E_{ij}^n \right.$$
$$\left. - \left(\frac{\Delta t_n}{\Delta t_{n-1}} \times \frac{2\Delta t_n + 3\Delta t_{n-1} + 3\Delta t_{n-2}}{\Delta t_{n-2}}\right)E_{ij}^{n-1} \right. \tag{42}$$
$$\left. + \left(\frac{\Delta t_n}{\Delta t_{n-2}} \times \frac{2\Delta t_n + 3\Delta t_{n-1}}{\Delta t_{n-1} + \Delta t_{n-2}}\right)E_{ij}^{n-2}\right]$$

For $U^*$ we have:

$$U_{ij}^{*\,n+1} = U_{ij}^{*\,n} + \frac{\Delta t_n}{6}\left[\left(\frac{\Delta t_n}{\Delta t_{n-1}} \times \frac{2\Delta t_n + 6\Delta t_{n-1} + 3\Delta t_{n-2}}{\Delta t_{n-1} + \Delta t_{n-2}} + 6\right)\left(F_{ij}^n + (F_{ij}^{*\,n})_t\right) \right.$$
$$\left. - \left(\frac{\Delta t_n}{\Delta t_{n-1}} \times \frac{2\Delta t_n + 3\Delta t_{n-1} + 3\Delta t_{n-2}}{\Delta t_{n-2}}\right)\left(F_{ij}^{n-1} + (F_{ij}^{*\,n-1})_t\right) \right. \tag{43}$$
$$\left. + \left(\frac{\Delta t_n}{\Delta t_{n-2}} \times \frac{2\Delta t_n + 3\Delta t_{n-1}}{\Delta t_{n-1} + \Delta t_{n-2}}\right)\left(F_{ij}^{n-2} + (F_{ij}^{*\,n-2})_t\right)\right]$$

To calculate $U^*$ we must also derive the equations for $F^*_t$. We will use the second order finite difference equations derived in the previous section:

$$(F_{ij}^{*\,n})_t = \frac{2\Delta t_{n-1} + \Delta t_{n-2}}{\Delta t_{n-1}(\Delta t_{n-1} + \Delta t_{n-2})}F_{ij}^{*\,n} - \frac{\Delta t_{n-1} + \Delta t_{n-2}}{\Delta t_{n-1}\Delta t_{n-2}}F_{ij}^{*\,n-1}$$
$$+ \frac{\Delta t_{n-1}}{\Delta t_{n-2}(\Delta t_{n-1} + \Delta t_{n-2})}F_{ij}^{*\,n-2} \tag{44}$$

$$(F_{ij}^{*\,n-1})_t = \frac{\Delta t_{n-2}}{\Delta t_{n-1}(\Delta t_{n-1} + \Delta t_{n-2})}F_{ij}^{*\,n} + \frac{\Delta t_{n-1} - \Delta t_{n-2}}{\Delta t_{n-1}\Delta t_{n-2}}F_{ij}^{*\,n-1}$$
$$- \frac{\Delta t_{n-1}}{\Delta t_{n-2}(\Delta t_{n-1} + \Delta t_{n-2})}F_{ij}^{*\,n-2} \tag{45}$$

$$(F_{ij}^{*\,n-2})_t = -\frac{\Delta t_{n-2}}{\Delta t_{n-1}(\Delta t_{n-1} + \Delta t_{n-2})}F_{ij}^{*\,n} + \frac{\Delta t_{n-1} + \Delta t_{n-2}}{\Delta t_{n-1}\Delta t_{n-2}}F_{ij}^{*\,n-1}$$
$$- \frac{\Delta t_{n-1} + 2\Delta t_{n-2}}{\Delta t_{n-2}(\Delta t_{n-1} + \Delta t_{n-2})}F_{ij}^{*\,n-2} \tag{46}$$

Plugging Eqs. (44) to (46) in (43) gives the equation to calculate $U^*$ in the next step. The equation for $V^*$ is derived similarly, and not represented here for the sake of brevity.

### 3.2.2.4 Implementation

We noticed that sudden increases in the time step sometimes lead to model instability. Therefore, we define a custom version of the exponential moving average and use it to set the value of $\Delta t_n$, as follows:

$$\Delta t_n := \begin{cases}\Delta t_n, & \Delta t_n \leq \Delta t_{n-1} \\ \alpha \Delta t_n + (1-\alpha)\Delta t_{n-1}, & \Delta t_n > \Delta t_{n-1}\end{cases} \tag{47}$$

where $α$ is a pre-defined coefficient between 0 and 1. We found $0.01 < α < 0.5$ to work well with Celeris. We call Eq. (47) *lazy* exponential moving average as it lets the time step to drop instantly, if required, but rise only gradually.

## 3.3  Boundary conditions

Two layers of ghost cells are considered at each boundary side and are used to implement the boundary conditions. Five types of boundary condition are implemented in Celeris Advent: fully reflective solid wall, sinewave maker [33], sponge layer [32], irregular wavemaker, and time series [36]. These boundary conditions can be applied to any of the four boundaries of the field.

## 3.4  Wave breaking

Wave breaking is not implemented in Celeris with a direct treatment. However, our experiments [33, 36] show that the numerical dissipation of the scheme caused primarily by using the minmod limiter imitates physical dissipation introduced by wave breaking. Kirby et al. [29] also discuss that in models with shock-capturing schemes, the implementation of an explicit formulation for breaking wave dissipation might be unnecessary. The MOST and GeoClaw models, commonly used in tsunami studies, are other examples in which numerical dissipation mimics wave breaking [2, 37, 39]. In the next sections we show that the wave breaking effect caused by the use of the limiter is able to adequately mimic the physical phenomenon as shown by comparisons of the numerical results with experimental measurements. As discussed before, our solver automatically reduces to the NLSW equations to continue simulating the runup on the beach.

## 3.5  Wet and dry cells

There is no definition for wet/dry cells in Celeris. All the cells are considered wet, though some with water depth of zero or close to zero. This treatment is possible thanks to the finite volume method and KP07 scheme and is favorable for our GPU implementation as it avoids branching in the GPU computations. To determine the runup or inundated area, one needs to define a minimum depth of inundation to distinguish wet and dry cells.

## 4.  Software development

Development of Celeris consisted of two major steps. Firstly, a robust solution for the mathematical model was developed such that it could be implemented on the GPU. The second step was the implementation of a user-friendly software which can drive the model. To fulfill the first step, we first implemented the derived equations and the mathematical model in MATLAB to validate them and improve them where necessary. We call this MATLAB series of the model Celeris Zero [32]. The purpose of developing Celeris Zero, was mostly validating the underlying mathematical model of the software, and not its implementation on the GPU. Development of the mathematical model and its implementation in Celeris Zero started in spring of 2014 and the first version was running in early 2015.

After we became confident that our model was suitable for our goal of developing an interactive and immersive coastal simulation software, we started the developments of the first official series of Celeris, called Celeris Advent [32], in winter 2015. Our goal in development of Celeris Advent was to provide a hassle-free software which can run on off-the-shelf Windows machines with minimum preparation. Therefore, we selected Microsoft's Direct3D library and its HLSL shader language to harness the power of the GPU, rather than general purpose GPU programming libraries such as CUDA, which requires some level of prior knowledge for preparation of the system and runs only on devices with NVIDIA GPU's. Celeris Advent is implemented mostly in C++ and HLSL, and it is an open-source code developed and redistributed under the terms of the GNU General Public License as published by the Free Software Foundation.

We released Celeris Advent (v.1.0.0) to the public in December 2016 [33]. At the time of writing this paper, Celeris Advent is downloaded over 2000 times by users from academia, industry, and government spanning over 50 countries and its user manual has been translated to Farsi, Spanish and Italian by independent users. Applications of Celeris Advent extend from research on coastal phenomena to recreational surf forecasts. One of our applications at University of Southern California, is a website which provides a five-day forecast of wave conditions at several US coasts, available at http://coastal.usc.edu/waves/. The adaptive time integration scheme is implemented in Celeris Advent (v.1.3.4) which is available to download at www.celeria.org. We recently introduced a new series of the Celeris software, called Celeris Base [34, 35], in which a modern game engine is used for implementation and immersive visualization capabilities are added to the software. Celeris Base is more suitable to researchers who would like to extend the capabilities of the software. We still recommend using Celeris Advent for general purpose simulations.

## 4.1  GPU implementation

Celeris harnesses the power of the GPU to run the TL17 numerical model as well as to concurrently visualize the results. We use shader programming for this purpose. The advantage of shaders compared to general purpose GPU programming languages such as CUDA is their portability between hardware and to some extent, between operating systems. Shaders are the core of 3D graphics rendering and game development. They were originated with the purpose of tweaking lighting effects in 3D rendering (hence the naming) but quickly became much more powerful and today are widely used in fixed-function rendering pipelines in graphics API's such as OpenGL and Direct3D.

The disadvantage of using shaders to implement our numerical solver is the need to restate the problem in terms of graphics primitives and setting up a dummy rendering pipeline. In Celeris Advent, we use Direct3D and its shader programming language, HLSL, to solve the governing equations on the GPU. This required setting up a dummy rendering pipeline to render a quad with two triangles and six vertices. We then divided the numerical scheme into smaller steps suitable to implement in pixel shaders. It was also required to fit the problem in the texture data structures. For example, the conservative variables vector in a cell, $\mathbf{U}_{ij}$=$[w_{ij}, P_{ij}, Q_{ij}]^{\mathrm{T}}$, are stored in *red*, *green*, *blue*, color components of a texel (i.e., a texture cell, similar to a pixel in a digital image). Therefore, a texture of size $(n_x+4)\times(n_y+4)$ on the GPU stores the conservative variables vectors for the solution domain of size $n_x \times n_y$, where four cells are added in each direction to account for the two layers of ghost cells on each boundary side. Other values in cells are also stored similarly in texture data structure, while globally constant variables are stored in constant buffers.

The explicit steps of the computation are relatively easy to implement on the GPU. That is because pixel shaders are designed to apply a kernel function on every texel of input textures and store the results in output textures. For example, the reconstruction step in KP07 is an explicit step, where the result depends only on some known values from the previous time step on the cell itself and its neighbors. But solving implicit equations such as Eq. (19) and Eq. (20) on the GPU can be challenging because the output on each cell is tied to the output of the other cells at the same time step and therefore the system needs to be solved simultaneously. We employed Cyclic Reduction (CR) algorithm to solve these equations on the GPU[33]. CR consists of two phases: forward reduction and backward substitution. In the forward reduction phase, the system is successively reduced to a smaller system with half the number of unknowns, until a system of 2 unknowns is achieved which can be solved trivially. In the backward substitution phase, the other half of the unknowns are found by substituting the previously found values into the equations. This process is well illustrated in Figure 4 of [33].

## 5. Model tests

### 5.1 Solitary wave run-up on a conical island

As the first validation test, we run a case which we previously used to validate Celeris Advent (v.1.0.0) with fixed time step in [33]. The experimental data of Briggs et al. [4] for solitary wave interaction around a conical island is frequently used to validate numerical models [11, 21, 38, 40] and has become a standard benchmark for Boussinesq-type models. The experimental setup consists of a circular island with 7.2 m base diameter and ¼ side slope, *s*, located in a 30m×25m wave tank with 0.32 m depth. Out of several cases in this set of experiments, we only test the case with target relative waveheights of $H/d$=0.20 which is expected to be more difficult for numerical models to simulate because of the higher non-linearity and the wave breaking condition. Briggs et al. [4] recorded the wave maximum run-up on the island and surface elevation time series on several gauges, which are used in this study for validation.

Our numerical setup for the conical island experiments consists of a 30m×30m domain with the conical island in the center and a solitary wave placed as an initial condition near the west boundary. The basin is extended by 5m to accommodate the solitary wave as an initial condition. The west and east boundaries are set to the sponge layer condition, while the north and south boundaries are fully reflective solid walls. The domain is discretized by 601×601 cells. We used adaptive time stepping with an initial time step of 0.0033s, corresponding to a CFL number of 0.145. No bottom friction is applied. The test case is performed with a slightly smaller relative waveheight at $H/d$=0.18. This reduction in waveheight ratio is used in several other studies such as [21], [40], and [11], as it is closer to the waveheight ratio observed downstream of the wavemaker.

Figure 5.1 shows the experimental setup and the gauges locations. Gauge #6 and #9 are in front of the island, while gauge #16 is on the side, and gauge #22 is behind the island. The numerical surface elevation compared to the experimental results are shown in Figure 5.2. The initial waveheight and subsequent draw-down are predicted well, which is consistent with numerical results from Celeris Advent with fixed time step [33] as well as results from other Boussinesq-type solvers [11, 21, 38, 40].

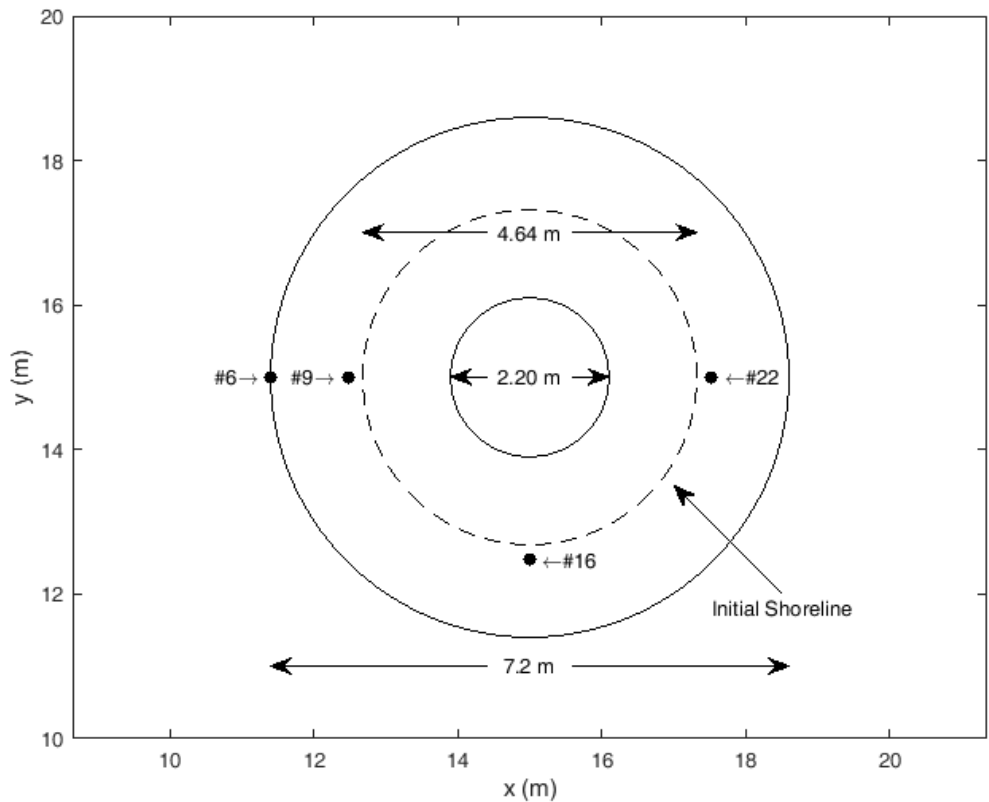

**Figure 5.1: Experimental setup of the conical island. The gauge locations are shown by dots and the wave approaches the island from the left.**

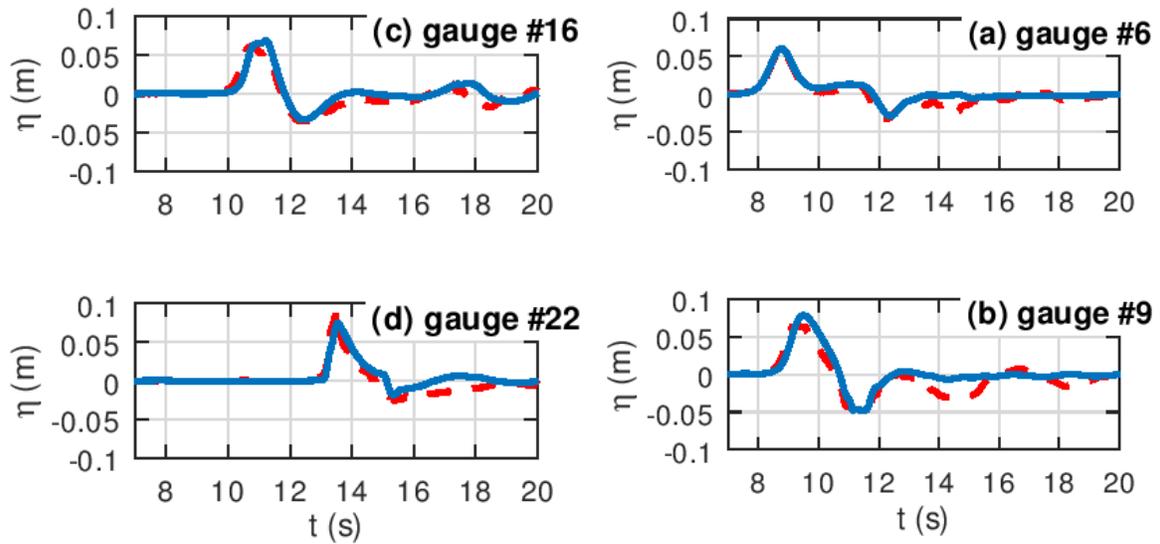

**Figure 5.2: Experimental (– –) and numerical (–) time series for Briggs et al.** [4] **benchmark at gauges #6, #9, #16, and #22 (a-d).**

Figure 5.3 compares the numerical and measured horizontal maximum run-ups on the island, scaled by the initial shoreline radius (2.32 m). Similar to [33] we used a threshold of $\delta = s\Delta x/3$ for water depth to determine the inundated area. The agreement between numerical data and measurements is very good even for the run-up on the back face of the island. As mentioned earlier, the wave breaks on the island and the strong agreement of data validates that numerical dissipation in TL17 successfully imitates wave breaking.

To compare our proposed adaptive scheme to the fixed time step scheme, we ran this experiment with $\Delta t = 0.0008s$ which was the largest fixed time step that resulted in a stable simulation. Figure 5.4 shows the adaptive time step variation over time and compares it to the fixed time step. As seen in this figure, adaptive time stepping let us run this experiment with a much larger time step during most of the simulation, saving a lot of computational effort.

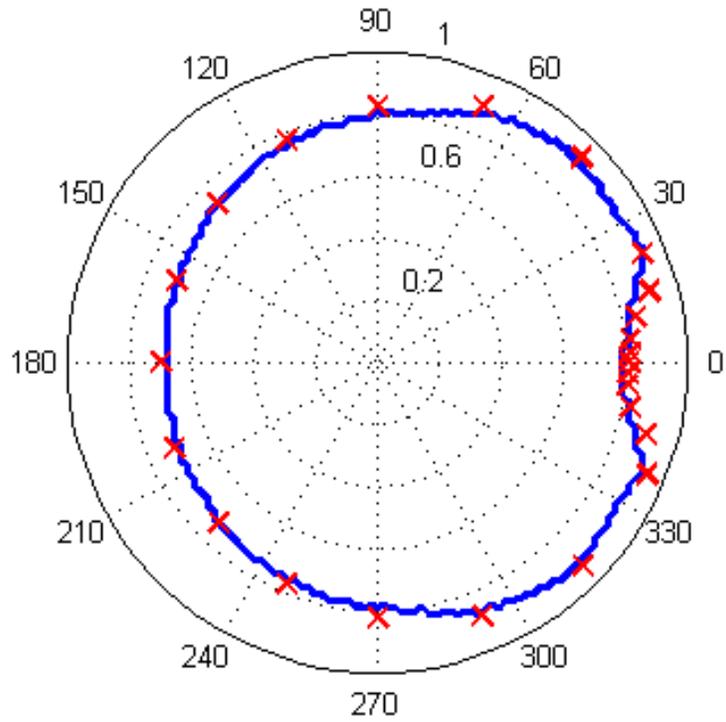

**Figure 5.3: Numerical (solid line) and measured (x) maximum horizontal run-up in Briggs et al. [4] benchmark.**

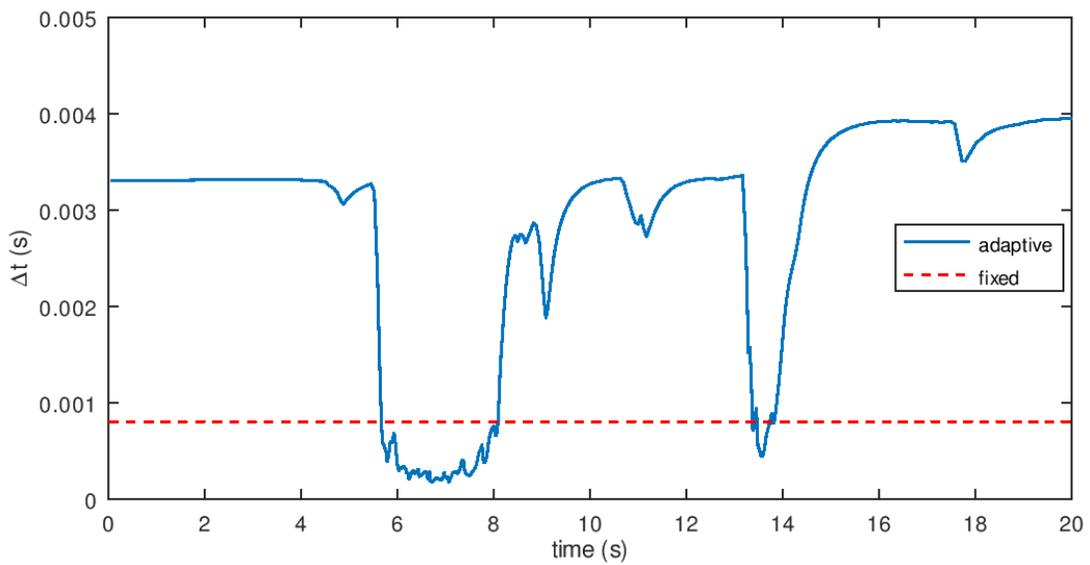

**Figure 5.4: Largest adaptive (solid line) and fixed (dashed line) time step for a stable simulation of Briggs et al. [4] benchmark.**

## 5.2 Breaking solitary wave runup on a slope with a conical island (Lynett et al., 2019)

We further test our proposed scheme by simulating the experiments of Lynett et al. [20] which have unfavorable hydrodynamic conditions for numerical models such as wave breaking and a moving shoreline. Such configuration is a reliable test for the robustness of our scheme. Lynett et al. [20] investigate the three-dimensional hydrodynamics of breaking solitary waves, traveling over shallow waters and their interaction with a conical island on a shelf. In a 26.5m×48.8m wave basin with 2.1m depth, a solitary wave is generated to propagate over a triangle-shaped shelf with a conical island as shown in Figure 5.5. This experiment is repeatedly simulated as a benchmark case for various Boussinesq-type models [9, 10, 29].

In our numerical experiment we constructed the bathymetry by combining the measured data of the shelf and the analytical equation of the conical island. The wave basin is discretized using a grid size of 0.1m×0.1m. The computational domain is further extended to the leeside of the wave maker, from $x$=0m to $x$=-10m, to accommodate the entire solitary wave as an initial condition. The waveheight is set to 0.39m at the wave generator location ($x$=0) where still water depth is set to 0.78m as in the laboratory experiment.

Figure 5.6 shows three snapshots of the simulated water surfaces at different times along with the corresponding local CFL numbers calculated from Eq. (22). The initial time step is set to $\Delta t$=0.0025s which later varies according to the evolving hydrodynamic conditions. The local CFL number at the corresponding time instants are shown in Figure 5.6. This figure shows that the local CFL number is entirely maintained below the critical value of 0.25 because of the adaptive time stepping. The free surface elevation simulated by Celeris Advent shows that the physical processes resulting from the wave evolution over the shallow shelf, such as wave steepening at the shelf front, wave scattering due to the island, wave runup, bore generation, wave breaking, and wave merging at the lee side of island, are well-captured in the simulation.

The time step variation during the simulation is shown in Figure 5.7 where the time instants corresponding to the snapshots in Figure 5.6 are shown by vertical dashed lines. The time step decreases when the solitary wave collides with the apex of the shelf and the island (i.e., Figure 5.6a). Then it reaches its minimum roughly at $t$=9s (i.e., Figure 5.6c and d) as the diffracted wave merges behind the island and thus the highest local Froude number is produced. Afterward, the wave condition becomes milder and time step is rapidly recovered up to 0.004s (i.e., Figure 5.6e and f). We also ran this experiment with a fixed time step and found $\Delta t$=0.0001 as the largest time step which still results in a stable simulation. This value is compared in Figure 5.7 with the time step of the adaptive scheme. As can be seen, the fixed time step requires a significantly smaller time step to keep the simulation stable, while the adaptive scheme only decreases the time step when the wave hits the island and is able to keep a relatively large time step for rest of the simulation.

Figure 5.8 compares calculated surface elevations with the laboratory data. For offshore-most four gauges (G1, G2, G4, G5), good agreement is seen between model and experimental data indicating that the model successfully predicts the wave transformation over three-dimensional shelf in front of the island. However, the more complicated hydrodynamic processes at G3, G6, G7, G8 and G9 due to the combined effects of wave breaking, bore formation, hydraulic jumps around the island and shoreline makes the predictions less reliable. Overall, the comparison between modelled and measured data shows reasonable agreement which is comparable to results from other studies [9, 10, 29].

Figure 5.9 compares the model results to the experiment for the flow velocity recorded at the measurement locations. For ADV1 which is deployed offshore of the island, the proposed scheme predicts $x$ direction velocity component, $u$, very well both in magnitude and phase, while keeping $y$ direction velocity component, $v$, close to zero consistent with the experiment. For ADV2, modelled $u$ and $v$ are generally consistent with the measurement, bearing in mind that some measurement data are not available due to the extremely dynamic wave motions. The comparison of the model and experiment at ADV3 shows good agreement, both for $u$ and $v$ as well.

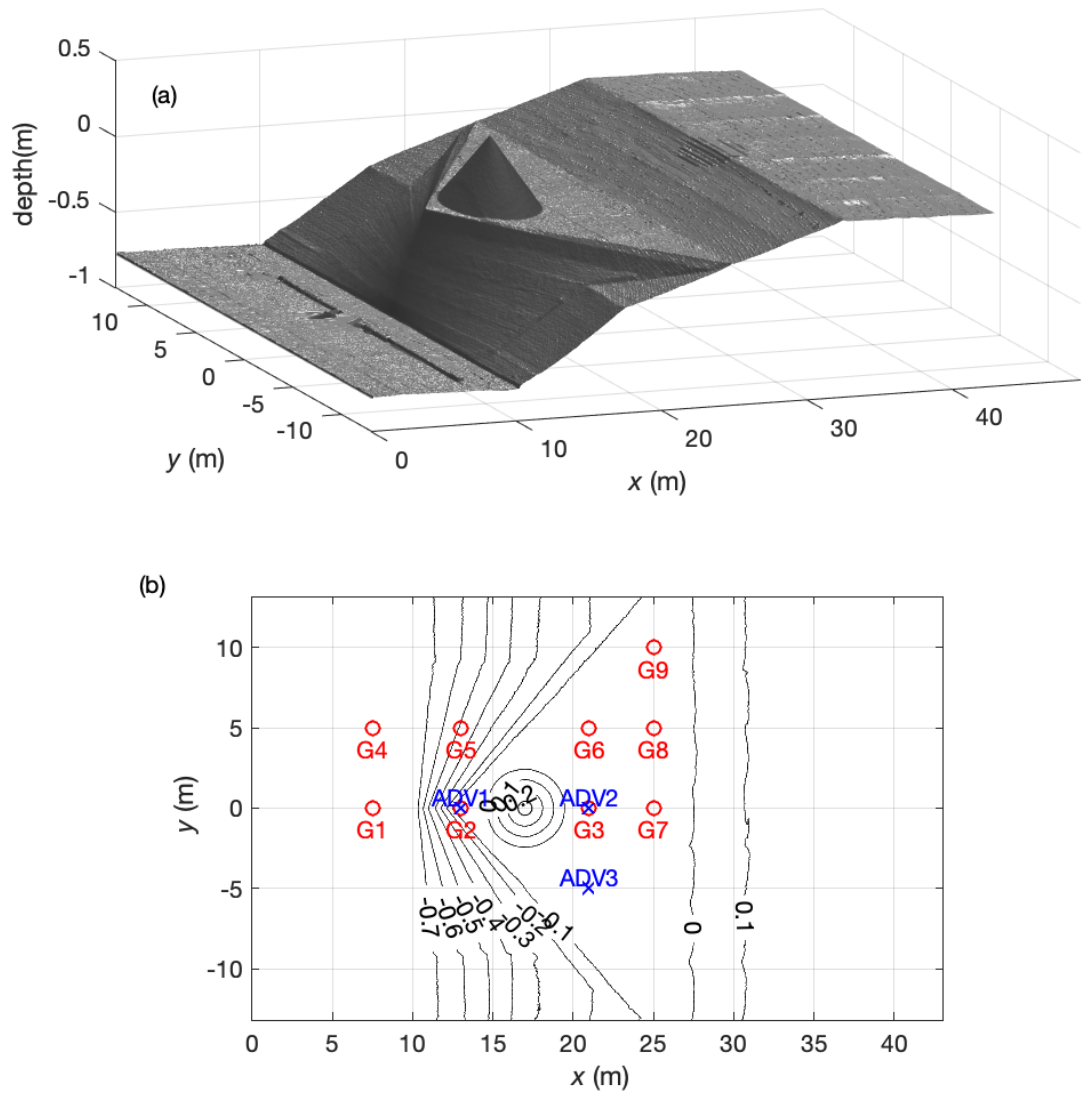

**Figure 5.5:** (a) Perspective view of bathymetry in the experiment of Lynett et al. (2019) and (b) bathymetric contours with apparatus setup. In a contour plot, red circles represent bottom pressure gauges while blue crosses indicate ADVs.

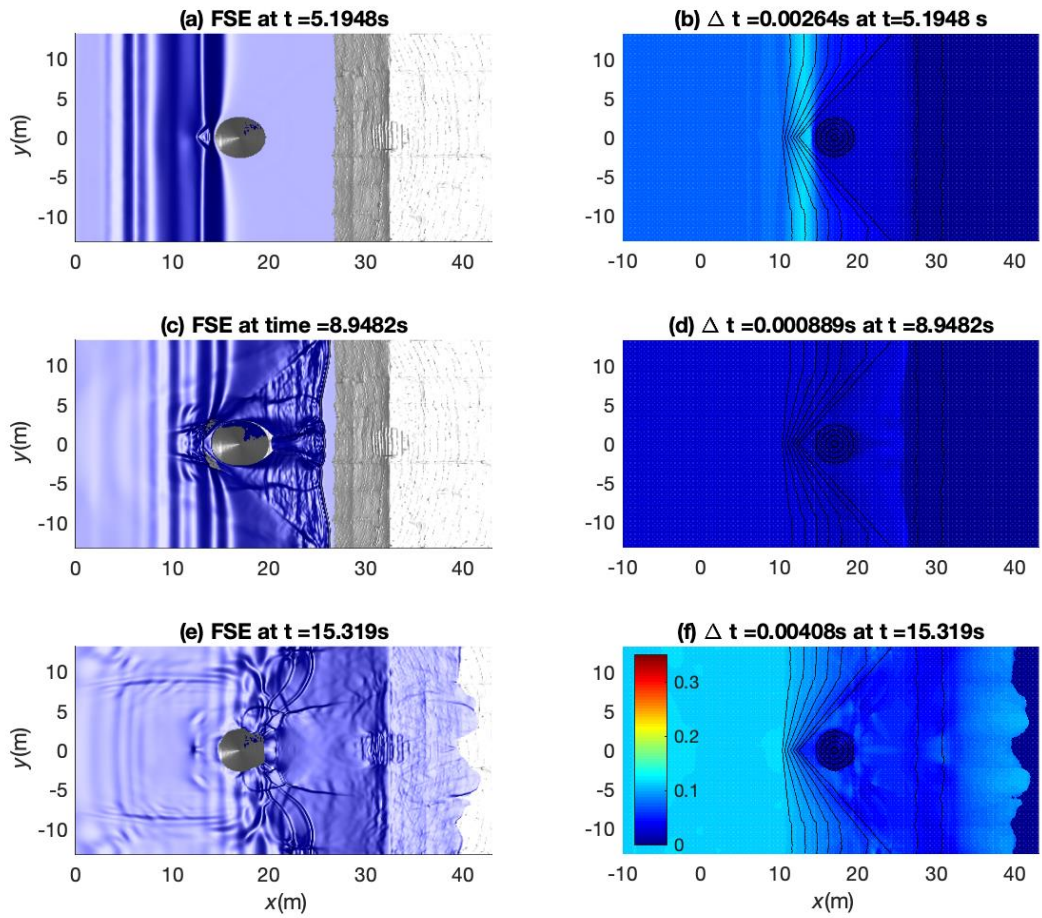

**Figure 5.6: Sequential snapshots of simulated free surface elevation (a, c, e) and corresponding map of local CFL number (b, d, f).**

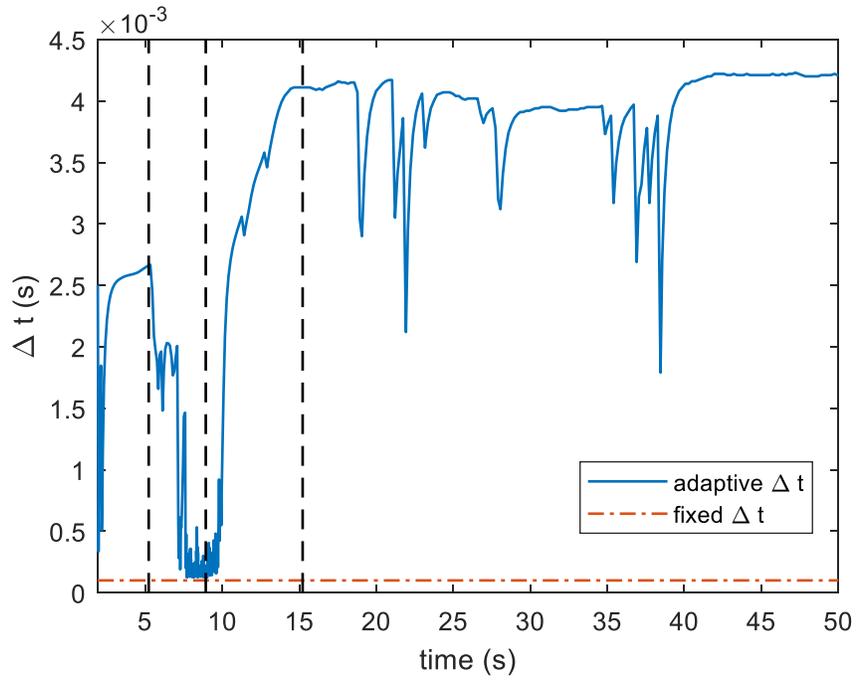

**Figure 5.7: Adaptive (solid line) and fixed (dashed-dot line) time steps during the simulation. Vertical dashed lines from left to right refer to the time sequences of Figure 5.6.**

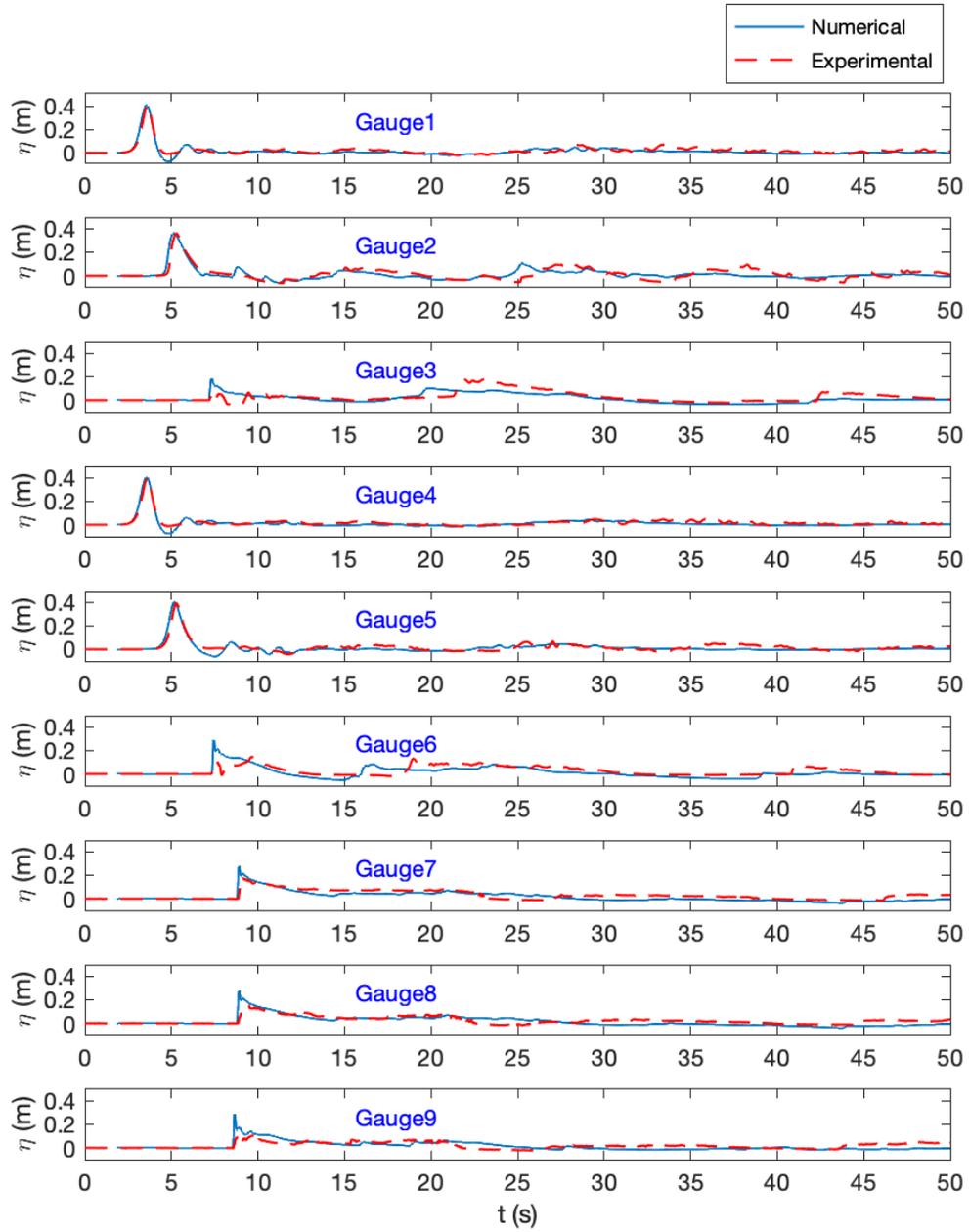

**Figure 5.8: Time series of free surface elevation at gauge locations.**

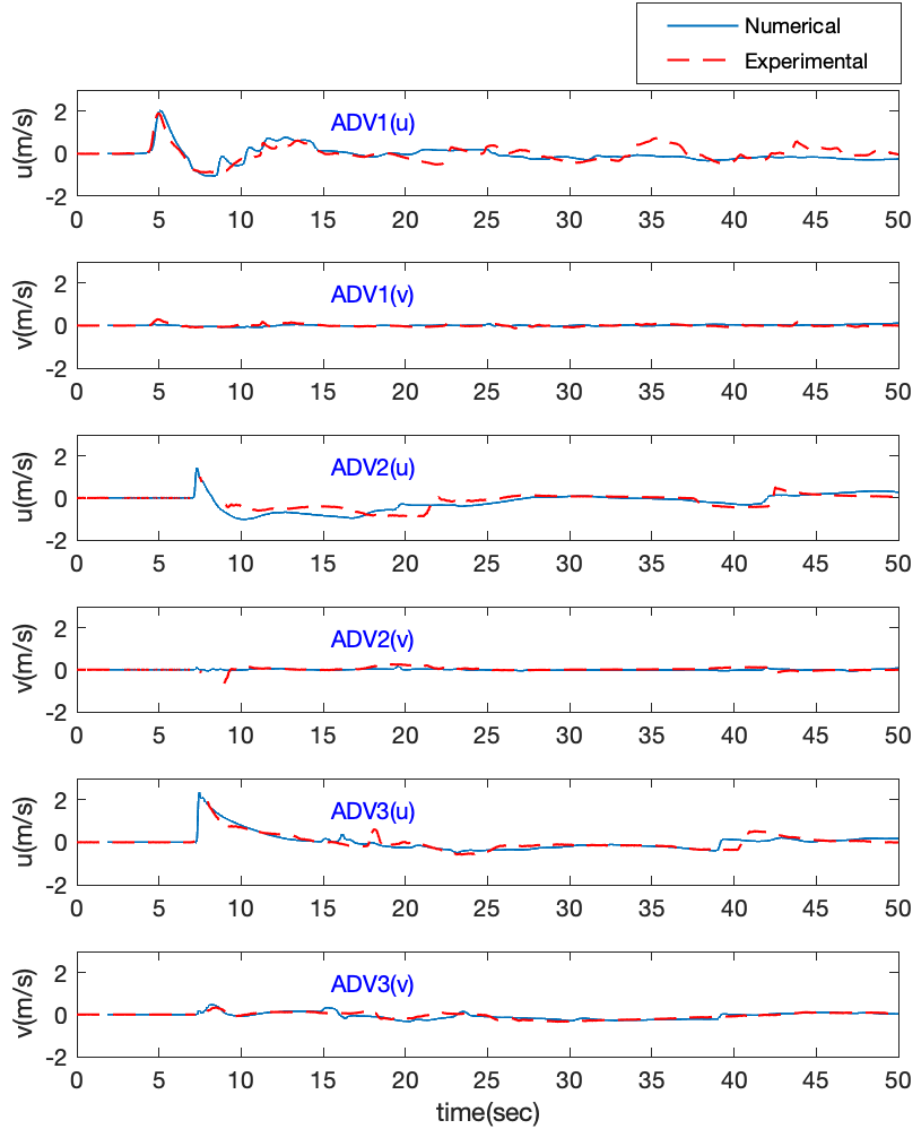

**Figure 5.9: Time series of velocity *u* and *v* at ADV location.**

## 5.3 Rip current with regular wave (Haller et al., 2002)

Another good benchmark for validating our proposed scheme is simulating a rip current experiment in which wave-current interactions are significant. It is widely known that Boussinesq-type models are capable of modelling wave-current interactions without incorporating short crested wave models [5, 6, 30, 31]. We simulate the laboratory experiments of Haller et al. [13] in this section.

Haller et al. [13] carried out laboratory experiments in a wave basin of size 18.2m×17m to investigate the rip current induced by wave breaking on the biplanar slope with a rip channel.

The bathymetry as used in the experiment is set up from the detailed dimensions of the experimental geometry while an analytical expression is borrowed for the submerged bar creation [1]. The reconstructed topography is shown in Figure 5.10. In the numerical experiment, the water depth is maintained at 0.678m at the wave-generation boundary and decreases until the water level sets shoreline at $x$=14.3m. A monochromatic wave with H=8.26cm and T=1.0s is generated on the left boundary and propagated onto the slope where two 7.32m-long alongshore submerged bars are located at $x$=12m with a spacing of 1.82m. Offshore current is anticipated to be generated through the rip channel due to the wave breaking induced momentum imbalance. A uniform grid size of 0.02m×0.02m is used for discretizing the domain while the adaptive time stepping with initial time step of 0.001s is applied. A bottom friction factor of 0.0025 is also considered for quadratic friction formula and the simulation time is set to 600s. For calculating time-averaged quantities of water level and velocity, simulated results are averaged over the last 500s which is equivalent to 500 wave periods. Figure 5.11 shows a snapshot of Celeris Advent simulating this experiment while visualizing vortical flows and wave breaking.

Figure 5.12 compares the alongshore variation of the calculated mean water level (MWL) at different locations with the measurements, showing good agreement. Figure 5.12a and Figure 5.12b indicate that wave setdown caused by increasing waveheights in the front of the bar system is well captured in the simulation and is consistent with the measurements. Figure 5.12c and Figure 5.12d show that the breaking induced wave setup taking place over the bar and setdown persisting through the rip channel are both well predicted by the model.

Figure 5.13 depicts the cross-shore variation of waveheight and MWL at the bar and at the rip channel. As pointed out by [12] the waveheight increases at the gap location due to the interaction between the incoming wave and offshore-directed rip current. Even though some discrepancy is noticed in the waveheight variation through the rip channel, the overall pattern of wave setup and setdown is well predicted by the model.

Figure 5.14 compares the time averaged cross-shore and longshore velocities at four different longshore transects with the laboratory data. The model results are generally consistent with the measurement as they well predict both cross-shore and longshore velocity variations over the rip channel system. Rip currents through the channel are clearly generated by the model as represented by offshore directed velocity (or negative value of $u_{\text{avg}}$) in Figure 5.14c and Figure 5.14d. The generated rip current also leads to longshore velocity changes which in turn form two opposite vortical flow pattern behind the rip channel [12, 13, 44]. These modelling results on rip current system including wave setup, wave setdown, and vorticity generation process confirm the capability of Celeris Advent and our proposed adaptive scheme.

Unlike previous experiments, the time step in this experiment did not significantly change during the simulation and only fluctuated within ±5% of its initial value. We attribute this small variation to the regular wave condition which was used in this experiment and the submerged bar which limited the run-up on the beach. Using adaptive scheme might not save us any computational time in experiments like this case, where extreme situations do not happen. However, the overhead of using the adaptive scheme is negligible and therefore we recommend always using this scheme over the fixed one.

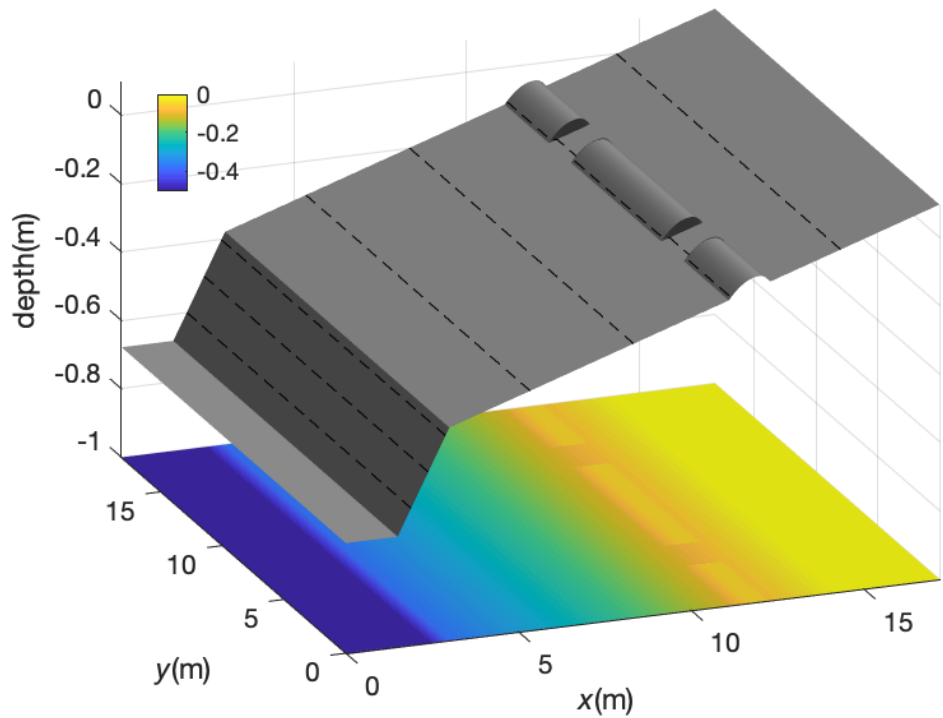

S

**Figure 5.10: Bathymetry of rip channel experiment in Haller et al. (2002). Dashed lines are contours at 0.1m depth intervals.**

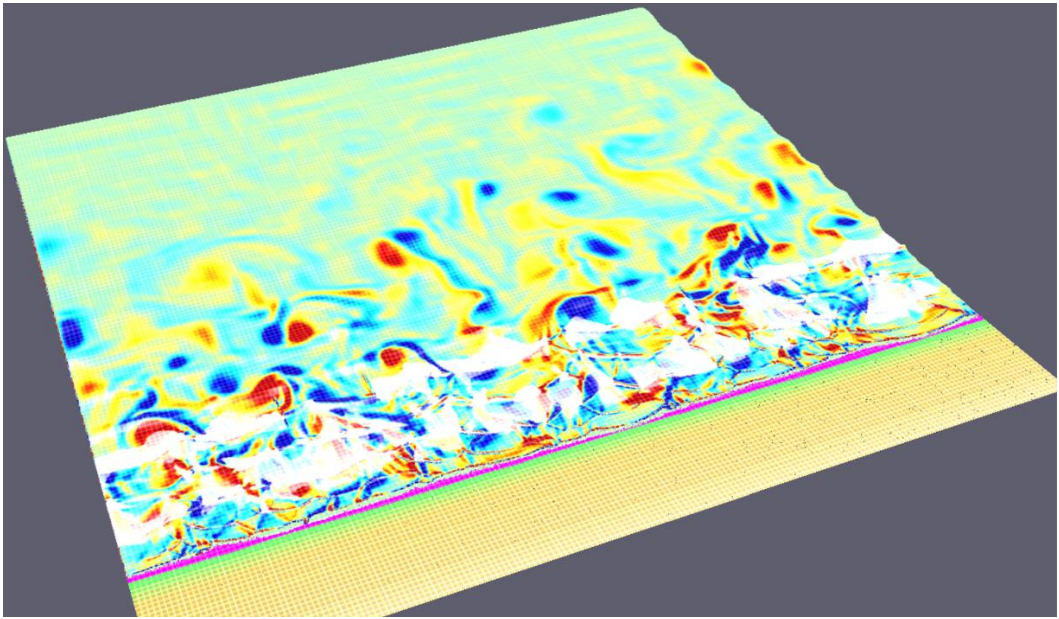

**Figure 5.11: Snapshot of Celeris Advent simulating the experiment of Haller et al. (2002) while visualizing the vorticities and wave breaking.**

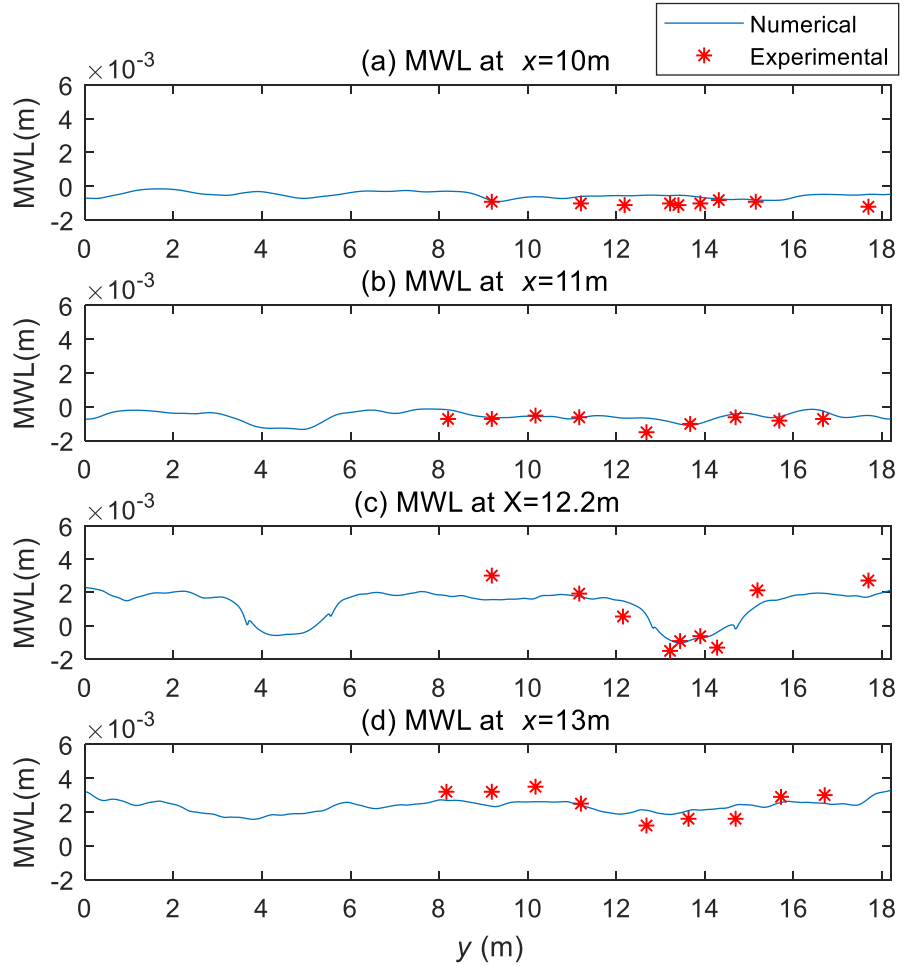

**Figure 5.12: Variation of mean water level (MWL) at four different longshore transects. Solid line represents calculated result and asterisk symbol denotes experimental data.**

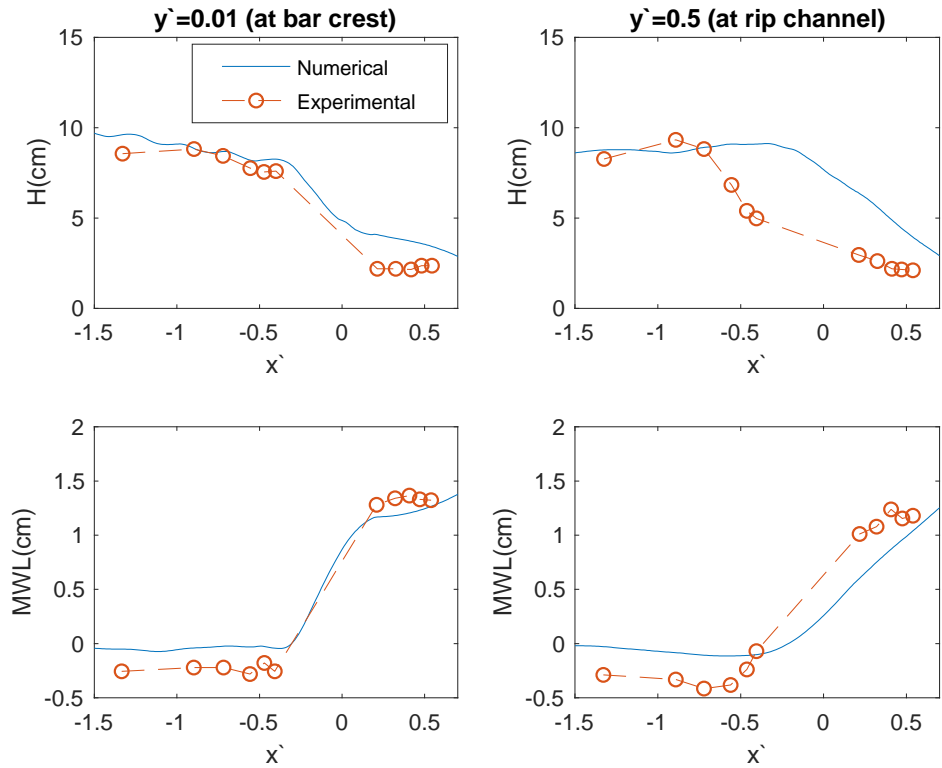

**Figure 5.13:** Variation of the waveheight and mean water level (MWL) at two different cross-shore transects. Solid line represents calculated result and dashed line with circle symbol denotes experimental data. Note that in the plots cross-shore distance *x* and alongshore distance y are nondimensionalized into x` and y` , respectively, according to [13].

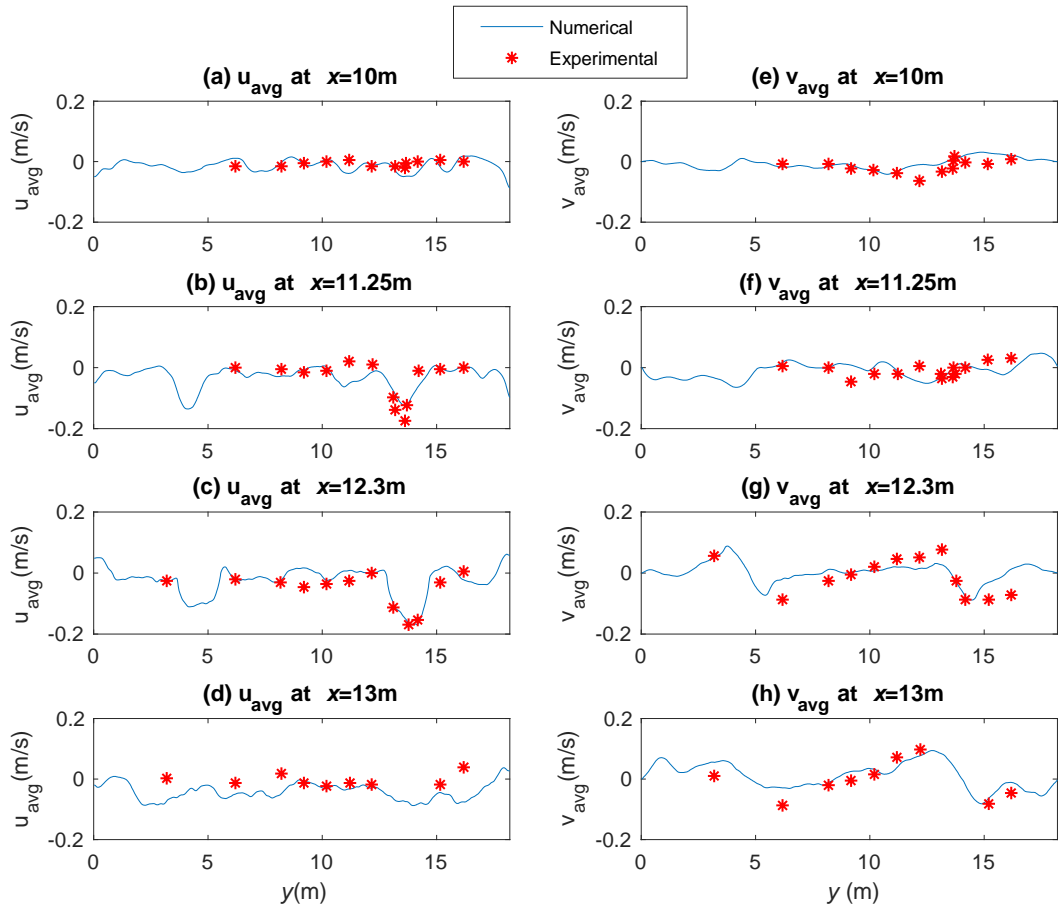

**Figure 5.14:** Variation of cross-shore ($u_{avg}$) and longshore ($v_{avg}$) time-averaged velocity component at four different longshore transects. Solid line represents calculated result, while asterisk symbol denotes experimental data. Panels (a-d) are for $u_{avg}$ and (e-h) for $v_{avg}$.

## 5.4 Rip current with irregular wave (Hamm, 1992)

The final case that we simulate is similar to the one in the previous section, but it includes irregular waves. We choose this test because the interaction of waves with a wide range of periods and the complex bathymetry is expected to result in challenging situations for a numerical model. For such situations, our adaptive scheme is more suitable for phase-resolving models since it enables them to control possible instabilities through spontaneous adjustment of the time step based on local Courant number [29]. This benchmark demonstrates how robust our proposed scheme is even when the wave condition is relatively harsh and random.

Hamm [14] experimentally investigated breaking-induced nearshore circulation in the wave basin where a plane beach with a rip channel was installed. Using both monochromatic and random waves, overall process of rip current generation by breaking waves was analyzed. The bathymetry in the experiment was created using the following equation.

$$z(x,y) = 0.1 - \frac{18-x}{30}\left[1 + 3\exp\left(-\frac{18-x}{3}\right)\cos^{10}\left(\frac{\pi y}{30}\right)\right] \tag{48}$$

where $z$ is the bottom elevation, $x$ is the onshore ward distance from the wave generator, $y$ is the longshore distance from the centerline of the channel. The bathymetry is shown in Figure 5.15. The domain is discretized by cells of size 0.02m×0.02m. The initial time step is set to 0.003s. Irregular waves are generated based on the JONSWAP spectrum with a significant waveheight of $H_s$=0.13m, and peak period of $T_p$=1.6s. The spectrum shown in Figure 5.16 is used to generate 68 wave components with discretized frequencies of $\Delta f$=0.01Hz. Quadratic formula is applied with a friction factor of 0.0025. Simulation time is set to 600s of which the beginning 100s is for spin up period and is not used in analyses.

Figure 5.17 shows the cross-shore variation of the significant waveheight, $H_s$, both at the plane beach and at the rip channel. It also shows the return current velocity, $u_{avg}$, at the rip channel. Calculated waveheights at two transects agree well with the measurement with some discrepancies, most possibly due to the imperfect reproducibility of the physical bathymetry [15]. Comparing Figure 5.17a and Figure 5.17b, the waveheight along the cross-shore transect at plane beach decreases more rapidly after $x$=11m than that along the rip channel, mostly due to the wave breaking taking place farther offshore along the plane beach transect. The returning current (i.e., rip current) induced by breakers at the rip channel is also well predicted by the model as shown in Figure 5.17c. The rip current interacts with the incoming waves in the opposite direction and consequently makes the waveheight larger around the channel as shown in Figure 5.17b.

Figure 5.18 shows the variation of the time step, $\Delta t$, during the simulation. It initiates from the starting value of 0.00325s and then steadily decreases as the wave propagates on to the shore. Unlike the regular wave case in which $\Delta t$ varies within a very limited bound and gets stabilized, $\Delta t$ of random wave simulation fluctuates continuously within a larger bound with occasional but drastic changes. This unpredictable variation is attributed to the random nature of irregular waves. The largest fixed time step that resulted in a stable simulation was $\Delta t = 0.0015$s which is almost half the average adaptive time step. This experiment did not have extreme events such as collision of a solitary wave on an island, yet the random nature of the waves required occasional drops in the time step to keep the simulation stable, which was possible thanks to the adaptive scheme.

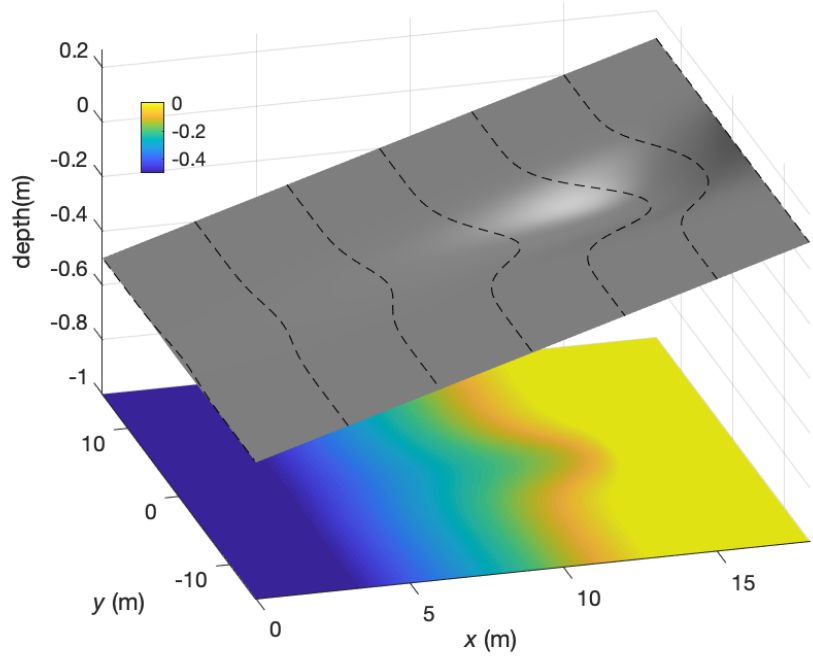

**Figure 5.15: Bathymetry of rip current experiment in Hamm (1993). Dashed lines are contours at 0.1m interval.**

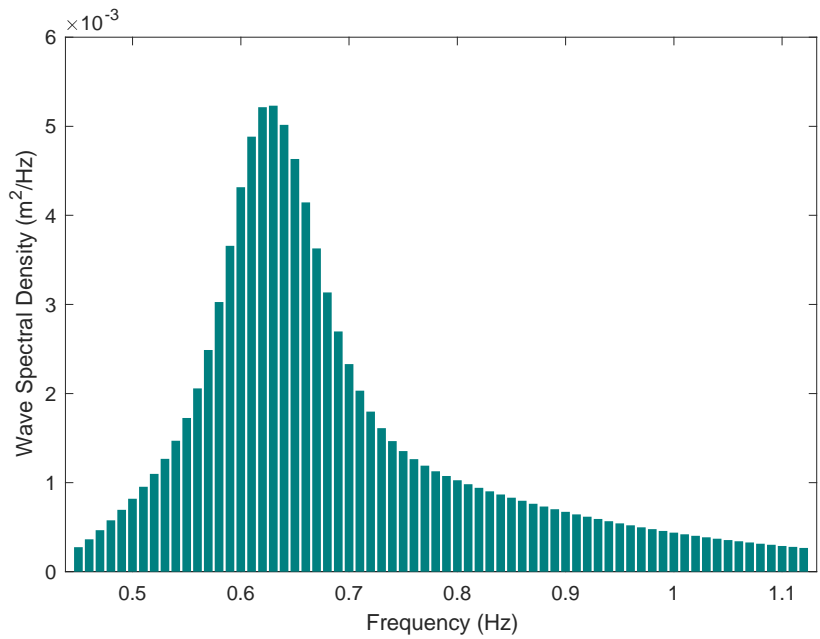

**Figure 5.16: JONSWAP wave spectrum discretized by 68 frequency components.**

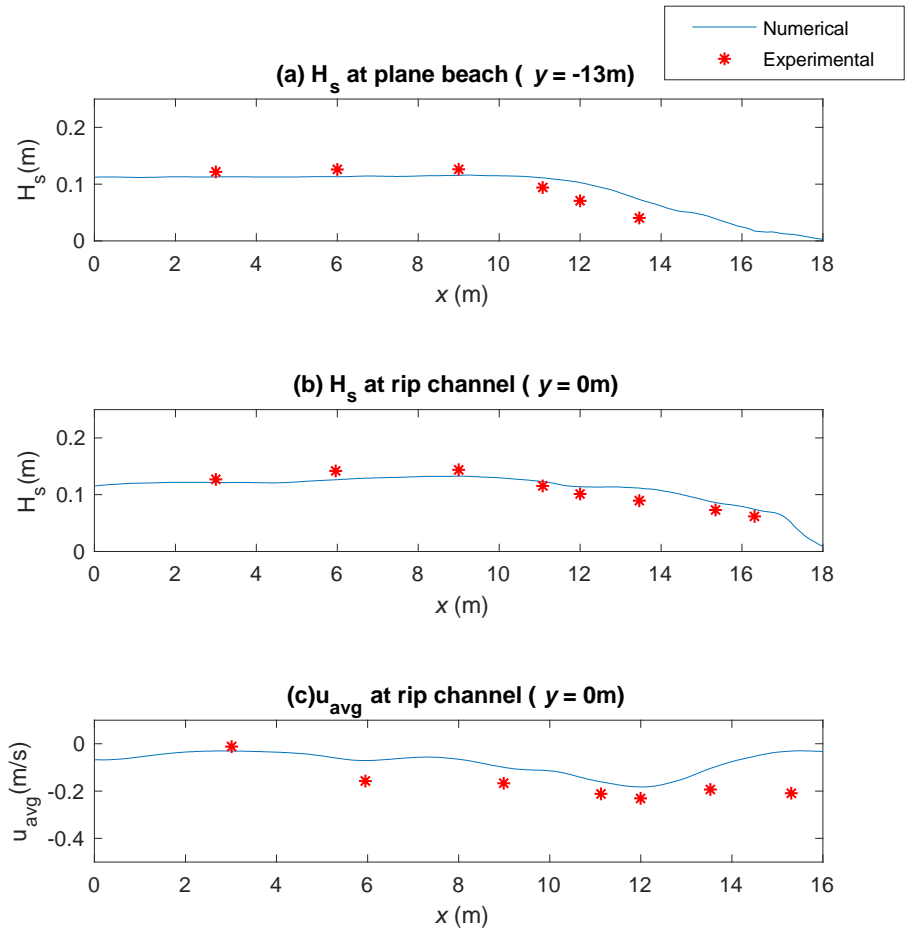

**Figure 5.17:** Cross-shore validation of the waveheight, Hs at (a) plane beach and (b) rip channel, and return flow, $u_{avg}$ at (c) the rip channel.

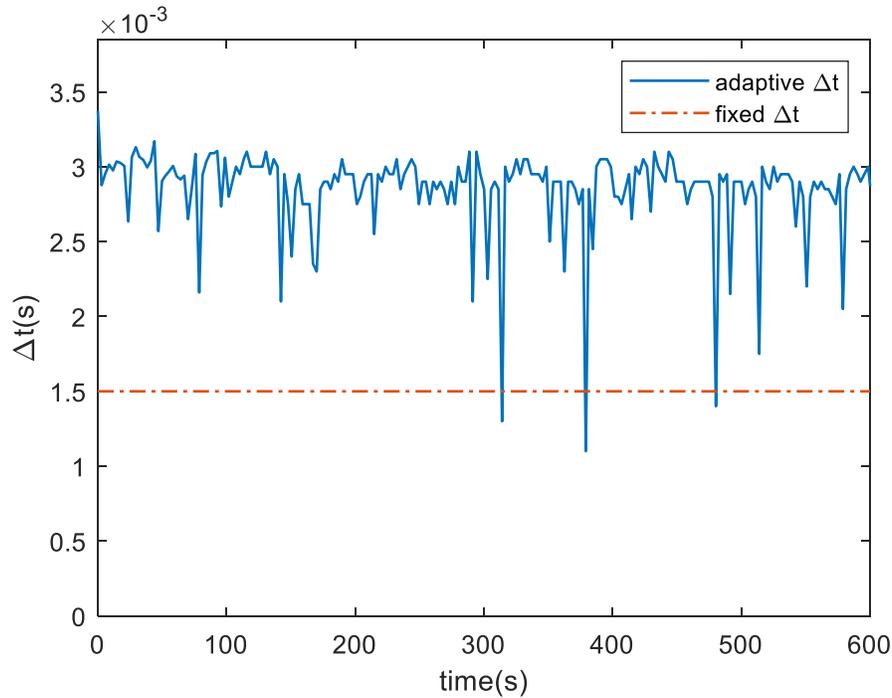

**Figure 5.18: Temporal variation of the time step. Solid and dashed-dot lines refer to adaptive and fixed time steps, respectively.**

# 6. Conclusion

We explained our development of the third order Adams Bashforth time stepping equations using the Lagrange formula for polynomial interpolation and assuming a variable time step value. We then employed these equations to solve the extended Boussinesq equations in time by developing the second order finite difference discretization equation for variable time steps and incorporating them in the rearrangement of the Boussinesq-type equations suggested by Wei and Kirby [42]. We implemented the resultant numerical scheme in the latest version of Celeris Advent (v1.3.4) and briefly explained this implementation. We validated the proposed adaptive scheme against several benchmarks proving the software's accuracy in modeling wave breaking, wave runup, irregular waves, and rip currents. The adaptive time stepping makes the model more robust by allowing it to keep the CFL number constant throughout the simulation. This is especially beneficial where the superposition of a wide range of wave conditions and a complex bathymetry (e.g., in field sites) creates occasional extreme conditions with large local Froude numbers. These brief extreme moments in the simulation are gracefully handled by the adaptive scheme using an accordingly small time step. As high-speed events diminish, the time step size is then recovered, letting the simulation continue efficiently.

# 7. Acknowledgements

This research was partially funded by the Office of Naval Research (ONR) award numbers N00014-13-1-0624 and N00014-17-1-2878.